\pdfoutput=1
\RequirePackage{ifpdf}
\ifpdf 
\documentclass[pdftex]{sigma}
\else
\documentclass{sigma}
\fi

\numberwithin{equation}{section}
\usepackage[all]{xy}

\newtheorem{thm}{Theorem}[section]

\theoremstyle{definition}
\newtheorem{defn}[thm]{Definition}

\newcommand{\bMor}{\bf Mor}
 \newcommand{\ba}{\begin{eqnarray}}
 \newcommand{\na}{\end{eqnarray}}
 
\newcommand{\Mor}{{\text{Mor}}}

 \newcommand{\IMor}{{\text{IMor}}}

 \renewcommand{\c}{\gamma}
 \newcommand{\e}{\varepsilon}

\def \mc{\mathcal}
\def \wt{\widetilde}

\def\E{\ifmmode{\mathbb E}\else{$\mathbb E$}\fi} %
\def\N{\ifmmode{\mathbb N}\else{$\mathbb N$}\fi} 
\def\R{\ifmmode{\mathbb R}\else{$\mathbb R$}\fi} 
\def\Q{\ifmmode{\mathbb Q}\else{$\mathbb Q$}\fi} 
\def\C{\ifmmode{\mathbb C}\else{$\mathbb C$}\fi} 
\def\H{\ifmmode{\mathbb H}\else{$\mathbb H$}\fi} 
\def\Z{\ifmmode{\mathbb Z}\else{$\mathbb Z$}\fi} 
\def\P{\ifmmode{\mathbb P}\else{$\mathbb P$}\fi} %
\def\T{\ifmmode{\mathbb T}\else{$\mathbb T$}\fi} %
\def\SS{\ifmmode{\mathbb S}\else{$\mathbb S$}\fi} %
\def\DD{\ifmmode{\mathbb D}\else{$\mathbb D$}\fi} %
\def\K{\ifmmode{\mathbb K}\else{$\mathbb K$}\fi}
\def\F{\ifmmode{\mathbb F}\else{$\mathbb F$}\fi} 

\begin{document}

\allowdisplaybreaks

\renewcommand{\thefootnote}{}

\newcommand{\arXivNumber}{2308.01595}

\renewcommand{\PaperNumber}{011}

\FirstPageHeading

\ShortArticleName{Twisted Sectors for Lagrangian Floer Theory on Symplectic Orbifolds}

\ArticleName{Twisted Sectors for Lagrangian Floer Theory\\ on Symplectic Orbifolds\footnote{This paper is a~contribution to the Special Issue on Differential Geometry Inspired by Mathematical Physics in honor of Jean-Pierre Bourguignon for his 75th birthday. The~full collection is available at \href{https://www.emis.de/journals/SIGMA/Bourguignon.html}{https://www.emis.de/journals/SIGMA/Bourguignon.html}}}

\Author{Bohui CHEN~$^{\rm a}$, Kaoru ONO~$^{\rm b}$ and Bai-Ling WANG~$^{\rm c}$}

\AuthorNameForHeading{B.~Chen, K.~Ono and B.-L.~Wang}

\Address{$^{\rm a)}$~School of Mathematics, Sichuan University, Chengdu 610064, P.R.~China}
\EmailD{\href{mailto:chenbohui@scu.edu.cn}{chenbohui@scu.edu.cn}}
\Address{$^{\rm b)}$~Research Institute for Mathematical Sciences, Kyoto University, Kyoto, 606-8502, Japan}
\EmailD{\href{mailto:ono@kurims.kyoto-u.ac.jp}{ono@kurims.kyoto-u.ac.jp}}

\Address{$^{\rm c)}$~Mathematical Sciences Institute, The National Australian University,\\
\hphantom{$^{\rm c)}$}~Canberra ACT 2601, Australia}
\EmailD{\href{mailto:bai-ling.wang@anu.edu.au}{bai-ling.wang@anu.edu.au}}

\ArticleDates{Received August 04, 2023, in final form January 16, 2024; Published online January 30, 2024}

\Abstract{The notion of twisted sectors play a crucial role in orbifold Gromov--Witten theory. We introduce the notion of dihedral twisted sectors in order to construct Lagrangian Floer theory on symplectic orbifolds and discuss related issues.}


\Classification{53D40; 53D37; 57R18}

\begin{flushright}
\begin{minipage}{80mm}
\it In honor of Professor Jean-Pierre Bourguignon\\ on the occasion of his 75th birthday
\end{minipage}
\end{flushright}

\renewcommand{\thefootnote}{\arabic{footnote}}
\setcounter{footnote}{0}

\section{Introduction}

The mathematical theory of quantum cohomology and Gromov--Witten theory on symplectic manifolds was first developed by Y.~Ruan and G.~Tian \cite{RT} in the semi-positive case.
This foundational work
was extended later by K.~Fukaya and the second named author~\cite{FO}, as well as J.~Li and G.~Tian~\cite{LT}, Y.~Ruan \cite{R} and B.~Siebert \cite{Si}.
W.~Chen and Y.~Ruan~\cite{CR} extended the theory further to the case of symplectic orbifolds, where the object called twisted sectors or inertia orbifolds plays a significant role.

The inception of the theory, now known as Floer theory, is attributed to A.~Floer. In \cite{Fl}, he laid the groundwork for the application to Lagrangian intersections, which became known as Lagrangian Floer theory.
After the study in the case of monotone symplectic manifolds due to Y.-G.~Oh~\cite{Oh,Oh-II}, K.~Fukaya, Y.-G.~Oh, H.~Ohta and the second named author~\cite{FOOO09-I,FOOO09-II}
constructed Lagrangian Floer theory for general Lagrangian submanifolds. In the case that Lagrangian submanifolds do not intersect the orbifold loci in a symplectic orbifold,
Lagrangian Floer theory was developed by C.-H.~Cho and M.~Poddar \cite{CP}.

However, when considering Lagrangians that may intersect the orbifold loci, complications arise.
To address this issue, we need to explore a specific class of Lagrangians and a variant of twisted sectors termed ``dihedral twisted sectors''.
These concepts are crucial for handling Lagrangian Floer theory. Since the diagonal in the product of two copies of an orbifold intersect the orbifold loci, we would like to treat the case
that Lagrangians may intersect the orbifold loci. The main objective of this article is to present and elucidate these important concepts, demonstrating their significance in the context of
Lagrangian Floer theory on symplectic orbifolds.\looseness=-1

The contents are as follows. In Section~\ref{sec2}, we briefly review Gromov--Witten theory and Lagrangian Floer theory on symplectic manifolds.
Moving on to Section~\ref{sec3}, we delve into fundamental notions concerning orbifolds and twisted sectors within the framework of orbifold Gromov--Witten theory.
In Section~\ref{sec4}, we introduce the notion of Lagrangians and their dihedral twisted sectors. In Section~\ref{sec5}, we give a brief discussion on the filtered $A_{\infty}$-algebra associated with a Lagrangian in a closed symplectic orbifold.

\section{The case of symplectic manifolds}\label{sec2}
In this section, we recall the framework of Gromov--Witten theory and Lagrangian Floer theory on symplectic manifolds.
We use the construction over the universal Novikov ring, which is defined by
\[\Lambda_0 = \Bigg\{ \sum_{i=1}^{\infty} a_i T^{\lambda_i} \mid a_i \in {\mathbb C}, \,\lambda_i \in {\mathbb R}_{\geq 0},\, \lambda_i \to \infty \  \text{for} \ i \to \infty\Bigg\}.\]
Its unique maximal ideal $\Lambda_+$ is defined by replacing the condition $\lambda_i \in {\mathbb R}_{\geq 0}$ by $\lambda_i > 0$.
The field of fractions of $\Lambda_0$ is denoted by $\Lambda$, which is called the universal Novikov field.
Let $(X, \omega)$ be a closed symplectic manifold and $J$ an almost complex structure compatible with $\omega$.
A~$J$-holo\-morphic map from a nodal Riemann surface $C$ equipped with distinct ordered $\ell$ marked points $\vec{x}=(x_1, \dots, x_\ell)$ away from nodes
to $(X,J)$ is called a stable map, if the automorphism group of $f\colon(C, \vec{x}) \to X$ is finite, i.e., the number of automorphisms $\varphi$ of $(C, \vec{x})$ with $f \circ \varphi = f$
is finite.

For $A \in {\rm H}_2(X;{\mathbb Z})$, we denote by ${\mathcal M}_{g, \ell}(X;A)$ the moduli space of stable maps of genus $g$,~$\ell$ marked points and representing $A$.
It is a compact metrizable space and carries the virtual fundamental class $[{\mathcal M}_{g, \ell}(X;A)]^{\rm vir}$.
Then the Gromov--Witten invariant is defined by
\[{\rm GW}_{g, \ell, A}\colon\ {\rm H}^*(X;{\mathbf Q})^{\otimes{^\ell}} \to {\mathbf Q}, (\alpha_1, \dots, \alpha_\ell) \mapsto \int_{[{\mathcal M}_{g, \ell}(X;A)]^{\rm vir}} {\rm ev}_1^* \alpha_1 \wedge \cdots\wedge
{\rm ev}_{\ell}^* \alpha_{\ell}\]
and it satisfies Kontsevich--Manin's axiom.\footnote{Except motivic axiom in symplectic case.}
In particular, using genus~0 Gromov--Witten invariant, one obtains quantum cup product $*_c$ parametrized by $c \in {\rm H}^*(X;\Lambda_0)$, i.e., the quantum cohomology ring of $(X, \omega)$.

For a Lagrangian submanifold $L$ in $(X, \omega)$, one can also consider $J$-holomorphic maps from bordered nodal Riemann surfaces $\Sigma$ with marked points to $(X,J)$ which map
the boundary $\partial \Sigma$ of $\Sigma$ to $L$. Here nodes and marked points are of two types, i.e., (1) on the interior of $\Sigma$, (2) on the boundary $\partial \Sigma$.
Boundary marked points $\vec{z}=(z_0, \dots, z_k)$ are disjoint from boundary nodes and interior marked points $\vec{x}=(x_1, \dots, x_\ell)$ are disjoint from interior nodes.
A holomorphic map~$u\colon(\Sigma, \partial \Sigma; \vec{z}, \vec{x}) \to (X,L)$ is called a bordered stable map, if the automorphism is finite.
For~$\beta \in {\rm H}_2(X,L;{\mathbb Z})$, denote by ${\mathcal M}_{k+1, \ell} (X,L;\beta)$ the moduli space of bordered stable maps of genus 0 and with $k+1$ boundary marked points,
$\ell$ interior marked points and connected boundary~$\partial \Sigma$, representing $\beta$.
For $k \geq 0$,\footnote{${\mathcal M}_{0, \ell} (X,L;\beta)$ may not be compact, unless we add stable map attached with a constant disk.} the moduli space ${\mathcal M}_{k+1, \ell} (X,L;\beta)$ is a compact metrizable space.
A spin structure of $L$, if exists, determines orientation on ${\mathcal M}_{k+1, \ell} (X,L;\beta)$ and one can construct a virtual fundamental chain.
Using the case that $\ell = 0$, we define a filtered $A_{\infty}$-algebra structure on the de Rham complex $\Omega^*(L) \hat{\otimes}_{\mathbb C} \Lambda_0$ by ${\mathfrak m}_k =\sum {\mathfrak m}_{k,\beta} T^{\int_\beta, \omega}$, where ${\mathfrak m}_{k,\beta}\colon \Omega(L)^k \to \Omega(L)$ is given by
\begin{gather*}
\begin{split}
&{\mathfrak m}_{k, \beta}(\xi_1, \dots, \xi_k) = (-1)^*({\rm ev}_0)_! ({\rm ev}_1^* \xi_1 \wedge \cdots \wedge {\rm ev}_k^* \xi_k), \qquad \text{unless} \quad (k,\beta)=(0, 0), (1, 0),\\
&{\mathfrak m}_{0,0}=0, \qquad{\mathfrak m}_{1, 0} ={\rm d} \qquad (\text{de Rham differential}).
\end{split}
\end{gather*}
Here $({\rm ev}_0)_!$ is the integration along fibers, if ${\rm ev}_0$ is a proper submersion. In general, it is defined using the theory of Kuranishi structures \cite{FOOO20}.

Using cases of all $\ell$, we defined operators $\mathfrak p$, $\mathfrak q$ (open-closed map, closed-open map) and bulk deformations of filtered $A_{\infty}$-structure by
a cycle in $X$ with $\Lambda_+$-coefficients.

For a cleanly intersecting pair $(L_0, L_1)$ of spin Lagrangian submanifolds, we can construct a~filtered $A_{\infty}$-bimodule over the filtered $A_{\infty}$-algebras associated to
$L_0$ and $L_1$.
We can extend these constructions for relatively spin Lagrangian submanifolds and relatively spin pair of Lagrangian submanifolds.
For the definition and discussion on
relative spin structures, see~\cite[Section~8.1.1]{FOOO09-II}.
The diagonal $\Delta_X \subset (X, -\omega) \times (X, \omega)$ is not necessary spin but relatively spin.
Set $\xi_1 \cup_Q \xi_2 = (-1)^{\deg \xi_1 (\deg \xi_2 +1)} {\mathfrak m}_2(\xi_1, \xi_2)$,
Then we have the following.

\begin{thm}[\cite{FOOOantisymp}] \label{ringiso} There is an isomorphism
\[I\colon\ ({\rm H}^*(X;\Lambda_0), *_0) \cong ({\rm H}^*(\Delta_X;\Lambda_0), \cup_Q).\]
\end{thm}

This statement must be plausible by naive comparison between the moduli spaces used for the product structures $*_0$, i.e., the small quantum product, and $\cup_Q$.
However, the stable map compactifications of these moduli spaces have different boundary structures. To rectify such a~discrepancy, a variant of the operation
$\mathfrak p$ mentioned above is used in the proof \cite{FOOOantisymp}.

\section{Orbifolds}\label{sec3}
\subsection{Presentation of an orbifold by a proper \'etale Lie groupoid}

A manifold is a geometric object locally modelled by a Euclidean space.
A geometric object locally modelled by a quotient space of a Euclidean space by a finite group action\footnote{We consider the action is effective in this article.} is
called an orbifold, which is a V-manifold\footnote{T.~Takakura asked Professor Satake what ``V'' stands for. His answer was Verzweigung.} introduced by Satake \cite{Sa}.
Namely, for each point $p$ on an $n$-dimensional orbifold $X$, there is a neighborhood $U$ of $p$ such that $U \cong B^n/\Gamma_p$, where
$B^n$ is an open ball $B^n \subset {\mathbb R}^n$ and $\Gamma_p$ is a finite group acting on ${\mathbb R}^n$ linearly.
We call $B^n \to B^n/\Gamma_p \cong U$ a local uniformizing cover.
In the case of manifolds, there are coordinate changes among local coordinate neighborhoods.
On an orbifold, there are equivariant diffeomorphisms among suitably shrunk local uniformization covers, which satisfy suitable compatibility condition.

We can define notions of (co)tangent vector bundles, more generally vector bundles, in a~natural way.
Differential forms and vector fields are defined as those on uniformization covers, which are invariant under the action of the finite groups~$\Gamma_p$.
For morphisms between orbifolds, if one defines it as a continuous map between the underlying topological spaces of orbifolds such that it is lifted to a smooth
equivariant map between uniformization covers, one may still pull back differential forms.
But this is not enough for pulling back vector bundles, in general. W.~Chen and Y.~Ruan introduced the notion of good maps.
Here we review the notion using the terminology of groupoids, e.g., \cite[Section~5.6]{MM}.

A groupoid is a category ${\mathcal C}=(C_0, C_1, s, t, m, u, i)$ such that all morphisms are invertible. Here~$C_0$ is the set of objects, $C_1$ is the set of morphisms,
$s \, (\text{resp. } t)\colon C_1 \to C_0$ is a map assigning the source (resp.\ target) to a morphism, \smash{$m\colon {C_1}_{s} \,\times_{C_0}\, {}_{t}{C_1} \to C_1$} is the composition
of morphisms.
Here \smash{${C_1}_{s} \,\times_{C_0}\, {}_{t}{C_1} = \{(f, g) \in C_1 \times C_1 \mid s(f)=t(g) \}$} and the composition $m$ enjoys the associativity.
From now on, we may write $m(f, g)= f \circ g$.
$i\colon C_1 \to C_1$ is a map assigning the inverse of a~morphism, $u\colon C_0 \to C_1$ is the map assigning the unit morphism to an object~in~$C_0$.

For $x, y \in C_0$, we define \smash{$x \,{\sim}_{\mathcal C}\, y$} if and only if $s^{-1}(x) \cap t^{-1}(y) \neq \varnothing$, namely, there is a morphism from $x$ to $y$.
Then \smash{${\sim}_{\mathcal C}$} is an equivalence relation.
Suppose that $C_0$, $C_1$ are manifolds, $s$, $t$ are smooth maps such that $s$ (or equivalently, $t$) is a submersion, \smash{${C_1}_{s}\, {\times}_{C_0}\, {}_{t}{C_1}$}
has also a structure of a manifold.
When $m$, $u$, $i$ are also smooth maps, $\mathcal C$ is called a Lie groupoid.

\begin{defn}\label{coarse} \quad
\begin{itemize}\itemsep=0pt
 \item[(1)] The quotient space of $C_0$ by the equivalence relation \smash{${\sim}_{\mathcal C}$} is called the coarse space of ${\mathcal C}$,
which is denoted by $\vert \mathcal C \vert$.
\item[(2)] A groupoid $\mathcal C$ is called a proper groupoid if $s \times t \colon C_1 \to C_0 \times C_0$ is proper.
\noindent
\item[(3)] A Lie groupoid $\mathcal C$ is called a \'etale Lie groupoid, if $s$ (or, equivalently, $t$) is a local diffeomorphism.
\end{itemize}
\end{defn}

Roughly speaking, for an orbifold $X$, we can construct a proper \'etale Lie groupoid ${\mathcal X}$ with an identification
$\vert {\mathcal X} \vert \cong X$ of the coarse space of $\mathcal X$ and the underlying topological space of the orbifold $X$.
Namely, $X_0$ is the disjoint union of local uniformization covers
and $X_1$ being the space of germs $[\psi]$ of local equivariant diffeomorphisms $\psi$ between suitable open subsets of local uniformization covers
such that $\psi$ induces the identity on an open subset of $\vert {\mathcal X} \vert \cong X$.
Denote by $\pi\colon X_0 \to \vert{\mathcal X} \vert$ the projection from the space of objects to the coarse space.

A proper \'etale Lie groupoid is locally described by an action groupoid (\!or translation groupoid\!) below \cite[Proposition 5.30]{MM}.

\begin{defn} Let $\Gamma$ be a group acting on $U$.
We set
$C_0 = U$, $C_1 = \Gamma \times U$, $s={\rm pr}_2$ (the projection to the second factor), $t\colon C_1=\Gamma \times U \to U$ the action of $\Gamma$ on $U$,
$u(x)=({\rm id}, x)$, $i(\gamma, x)=\bigl(\gamma^{-1}, \gamma \cdot x\bigr)$.
Define $m((\gamma, x), (\sigma, y))=(\gamma \cdot \sigma, y)$, when $x=\sigma \cdot y$.
Then ${\mathcal C}=(C_0, C_1, s, t, m, u, i)$ is a groupoid, which is called an action groupoid and is denoted by $\Gamma \ltimes U$.
When $U$ is a manifold and a finite group $\Gamma$ acts on $U$ smoothly, $\mathcal C$ is a proper \'etale Lie groupoid.
\end{defn}

From now on, we use the presentation of an orbifold $X$ by a proper \'etale Lie groupoid ${\mathcal X}$ in the following arguments.

\begin{defn}\label{def3.3} \quad
\begin{itemize}\itemsep=0pt
\item[(1)] A differential form on an orbifold $X$ presented by a proper \'etale Lie groupoid ${\mathcal X}$ is a~pair of differential forms
$\eta_0$, $\eta_1$ on $X_0$ and $X_1$, respectively, such that $\eta_1=s^*\eta_0=t^*\eta_0$. In particular, a symplectic form on $X$ is a pair of symplectic forms
$\omega_0$, $\omega_1$ on $X_0$, $X_1$,
respectively, such that $\omega_1=s^*\omega_0=t^*\omega_0$.
\item[(2)] An almost complex structure on orbifold $X$ is a pair of almost complex structures $J_0$, $J_1$ on $X_0$, $X_1$, respectively, such that
$s_* \circ J_1 = J_0 \circ s_*, t_* \circ J_1 = J_0 \circ t_*$.
\item[(3)] A vector bundle $E$ on an orbifold $X$ is a pair of vector bundles $E_0$, $E_1$ on $X_0$, $X_1$, respectively, equipped with consistent isomorphisms\footnote{These isomorphisms gives an action of $X_1$ on $E_0$.}
$s^*E_0 \cong E_1$, $t^*E_0 \cong E_1$.
Principal bundles on an orbifold is defined in the same manner.
\end{itemize}
\end{defn}

Next, we discuss the notion of morphisms between orbifolds.
Let $\mathcal X$ and $\mathcal Y$ be proper \'etale Lie groupoids representing orbifolds $X$ and $Y$.
Since a groupoid is a category, it is natural to consider a functor $F=(F_0, F_1)$ such that $F_i\colon X_i \to Y_i, \ i=0,1$ are smooth (smooth functor).
 We call such a functor a strict smooth morphism from $\mathcal X $ to $\mathcal Y$.
It is, however, not sufficient, since, even in the case of smooth maps between manifolds, the image of a coordinate chart of $X$ is not necessarily contained in
a coordinate chart in $Y$. Therefore, we need to take a refinement of a groupoid.

\begin{defn}
A refinement of a proper \'etale Lie groupoid ${\mathcal X}$ associated with an open covering $\bigl\{ U^{(j)} \bigr\}$ of $X_0$ is a proper Lie groupoid
with $U_0 = \bigsqcup_j U^{(j)}$, $U_1 = \bigsqcup_{i,j} t^{-1}\bigl(U^{(j)}\bigr) \cap s^{-1}\bigl(U^{i}\bigr)$ such that the structure maps $s$, $t$, $m$, $u$, $i$ are naturally induced from those
for ${\mathcal X}$.
\end{defn}

For orbifolds $X$, $Y$, we define a morphism from $X$ to $Y$ as a smooth functor from some refinement of $\mathcal X$ to $\mathcal Y$.
Let us go back to the case of smooth maps between manifolds, the {\it same} map may be described as various system of maps between
coordinate charts, i.e., the image of a~coordinate chart of $X$ may be contained in various coordinate charts of $Y$.
Hence the description as a smooth functor is not unique. Thus we need to consider smooth natural transformation between two smooth functors.

To be precise, let $\bMor_0({\mathcal X}, {\mathcal Y})$ be the object space consisting of smooth functors from a~refinement of ${\mathcal X}$ to $\mathcal Y$,
with its element given by
\[\xymatrix{
 {\mathcal X} & \ar[l]_{\phi } {\mathcal U} \ar[r]^u &{\mathcal Y}, }
\]
where $\phi\colon {\mathcal U} \to {\mathcal X}$ is a refinement of ${\mathcal X}$, and $u\colon {\mathcal U} \to {\mathcal Y}$ is a strict smooth morphism. We simply denote this object by $({\mathcal U}, \phi, u)$.
Given two objects $({\mathcal U}, \phi, u)$ and $( {\mathcal V}, \psi, v)$ a morphism from $({\mathcal U}, \phi, u)$ and $( {\mathcal V}, \psi, v)$ is a common refinement ${\mathcal W}$ of
${\mathcal U}$ and ${\mathcal V}$ together with a natural transformation
$ \alpha\colon u\circ \pi_1 \Longrightarrow v\circ \pi_2 $
as illustrated in following diagram:
\[
\xymatrix{
 &&{\mathcal U} \ar@/_/[lld]_{\phi} \ar@/^/[rrd]^{u} \ar@2[dd]_{\alpha} && \\ {\mathcal X} & \ar[l] {\mathcal W} \ar[ru]^{\pi_1} \ar[rd]_{\pi_2} &&& {\mathcal Y} . \\
 && \ar@/^/[llu]^\psi {\mathcal V} \ar@/_/[rru]_{v}& &
}
\]
This forms the morphism space $\bMor_1(({\mathcal U}, \phi, u), ( {\mathcal V}, \psi, v)) \subset \bMor_1({\mathcal X}, {\mathcal Y})$.
The composition of two composable morphisms and other structure maps can be found in~\cite{CDL} where the Sobolev completion of $\bMor({\mathcal X}, {\mathcal Y})$ is developed.\footnote{There is a formulation in terms of
 bibundles. For the purpose of
the moduli spaces of stable maps with Kuranishi structures, we use the description given here.}
For any continuous morphism $f \in \bMor_0({\mathcal X}, {\mathcal Y})$ in this sense, we can pull-back vector bundle $E$ in the sense of Definition \ref{def3.3}\,(3).
For $f, g \in \bMor_0({\mathcal X}, {\mathcal Y})$, if there is a morphism $T \in \bMor_1({\mathcal X}, {\mathcal Y})$ from $f$ to $g$, $T$ induces an isomorphism between~$f^*E$ and~$g^*E$.

\subsection{Twisted sector and orbifold stable maps}
We introduce the notion of the twisted sector (or inertia groupoid) for an orbifold $X$ presented by a proper \'etale Lie groupoid $\mathcal X$.

\begin{defn} For a proper \'etale Lie groupoid $\mathcal X$, we set
\begin{gather*}
IX_0=\{ a \in X_1 \mid s(a) = t(a)\},\qquad
IX_1=\bigl\{a \overset{g}{\to} b \mid a, b \in IX_0,\, g \in X_1, \, b =g \circ a \circ i(g) \bigr\},\\
s \bigl(a \overset{g}{\to} b\bigr) = a, \qquad t \bigl(a \overset{g}{\to} b\bigr) = b,\qquad  m\bigl(b \overset{h}{\to} c, a \overset{g}{\to} b\bigr) = a \overset{h \circ g}{\to} c,\\
i\bigl(a \overset{g}{\to} b\bigr)=b \overset{i(g)}{\to} a,\qquad u(a)=a \overset{u(s(a))}{\to} a.
\end{gather*}
Then we obtain a groupoid
\[\mathcal{IX}=(IX_0, IX_1, s, t, m, u, i).\]
We call it the twisted sector (or inertia groupoid) of $\mathcal X$. The twisted sector of a proper \'etale Lie groupoid $\mathcal X$
is a proper \'etale Lie groupoid, although the dimension depends on its connected components.
Note that $IX_0$ contains the space of identities, which is identified with $X_0$.
The restriction of ${\mathcal IX}$ to $X_0$ is called the trivial twisted sector (or the untwisted sector).
\end{defn}

For example, the twisted sector ${\mathcal I}(\Gamma \ltimes U)$ of an action groupoid $\Gamma \ltimes U$ is described as
\begin{gather*}
I(\Gamma \ltimes U)_0=\{(\gamma, x) \in \Gamma \times U \mid \gamma \cdot x = x\}, \\
I(\Gamma \ltimes U)_1= \bigl\{(\gamma, x) \overset{\rho}{\to} (\sigma, y) \mid \gamma \cdot x = x,\, \sigma \cdot y = y,\, y=\rho \cdot x, \,\sigma = \rho \circ \gamma \circ \rho^{-1} \bigr\}.
\end{gather*}

When an orbifold $X$ is equipped with an almost complex structure $J=(J_0, J_1)$, we assign the age (or the degree shifting number) to $a \in IX_0$.
Considering the situation locally, we reduce the discussion to the case of an action groupoid.
Then $a$ is given by $(\gamma, x) \in I(\Gamma \ltimes U)_0$. We regard $\gamma\colon U\to U$ around $x$ a $J_0$-linear map $\gamma\colon {\mathbb C}^n \to {\mathbb C}^n$
($x$ corresponds to the origin of $ {\mathbb C}^n$).
Since $\Gamma$ is a finite group, there is a minimal positive integer $m$ such that $\gamma^m = \rm{id}$ and the action by $\gamma$ around $x$ is diagonalizable and
presented by a matrix conjugate to
\[\operatorname{diag}\bigl(\exp\bigl(2 \pi m_1 \sqrt{-1}/m\bigr), \dots, \exp\bigl(2 \pi m_n \sqrt{-1}/m\bigr)\bigr),\]
where $0 \leq m_j < m$, $j=1, \dots, n$.
Then we set ${\rm age}(\gamma, x)= (m_1+ \dots + m_n)/m$ and call it the age (or the degree shifting number) of $(\gamma, x)$.
The age is a locally constant function on (the space of objects of) the twisted sector.
More generally, for an orbifold complex vector bundle, we define the age of the $\Gamma$-action on the fiber.

The age is an important invariant associated with the action of $\gamma \in \Gamma$ on a complex vector bundle.
Specifically, when considering the Dolbeault operator acting on sections of holomorphic orbifold vector bundle over an orbifold Riemann surface,
its Fredholm index is, by definition, an integer.
However, the orbifold Chern number is a rational number, not necessarily an integer.
Consequently, in this context, the topological quantity appearing in Riemann--Roch theorem is corrected by the age (for the local action on the vector bundle).
In the theory of stable maps from an orbifold Riemann surface, an elliptic operator acting on the pull-back of the tangent bundle of the target by
an orbifold stable map appears as the linearization of the equation for pseudo-holomorphic maps. Such an operator has the same principal symbol
as the Dolbeault operator with coefficients in a holomorphic orbifold vector bundle, i.e., the pull-back of the tangent bundle.
Thus the (virtual) dimension of the moduli space is expressed by the Fredholm index of the linearization operator and the age mentioned above plays an important role.

A holomorphic map from an orbifold Riemann surface $C$ to an almost complex orbifold~$X$ is defined to be a smooth functor $F=(F_0, F_1)$ from a refinement $\mathcal C$ of
a proper \'etale Lie groupoid representing $C$ to a proper \'etale Lie groupoid $\mathcal X$ equipped with an almost complex structure~$(J_0, J_1)$ representing $X$ such that
$F_0$, $F_1$ are holomorphic with respect to $J_0$, $J_1$, respectively.
Here, $C$ is an effective orbifold (locally, the action groupoid $\Gamma \ltimes D$ associated with an effective action of a finite group $\Gamma$
on the unit disk $D \subset \mathbb C$
keeping the origin fixed). The almost complex orbifold $X$ is also locally presented by an action groupoid associated with an effective action of a finite group $G$ on
an almost complex manifold. Then a morphism $\phi\colon \Gamma \ltimes D \to G \ltimes U$ is given by a pair of $\phi_0\colon D \to U$ and a homomorphism $\phi_1\colon \Gamma \to G$
such that $\phi_0$ is $\phi_1$-equivariant $J_0$-holomorphic map. In orbifold Gromov--Witten theory, we only consider those such that the homomorphism
$\phi_1$ is injective.

We next explain a pseudo-holomorphic $\phi\colon {\mathcal C} \to {\mathcal X}$ takes values in the twisted sector.
Pick an action groupoid $\Gamma \ltimes D$ (with $O \in D$ representing $p$) presenting a neighborhood of $p \in C$ and
an action groupoid $G \ltimes U$ presenting a neighborhood of $\vert \phi \vert (p) \in X$.
 Here $\vert \phi \vert$ is the map from~$\vert \mathcal C \vert$ to $\vert \mathcal X \vert$ induced by $\phi$.
 Let $\eta \in \Gamma$ be the generator corresponding to ${\rm exp}\bigl(2\pi \sqrt{-1}/m\bigr)$. Then~$(\phi_1(\eta), \phi_0(O))$ belongs to $IX_0$.
If $\eta$ is not trivial, $(\phi_1(\eta), \phi_0(O))$ belongs to a non-trivial twisted sector, since $\phi_1$ is injective.
When there is a natural transformation between $\phi=(\phi_0, \phi_1)$ and $\phi'=(\phi_0', \phi_1')$, the equivalence classes of $(\phi_1(\eta), \phi_0(O))$ and $(\phi_1'(\eta), \phi'_0(O))$
in $\vert I{\mathcal X} \vert$
are the same. Therefore, the image in $\vert {\mathcal IX} \vert$ is well defined.

An orbifold stable map from an orbifold nodal Riemann surface $C$ to an almost complex orbifold $X$ is defined in the following way.
Take a normalization $p\colon \widetilde{C} \to C$ of $C$. Then, for each (possibly orbifold) node $z \in C$, there is a pair $(\tilde{z}', \tilde{z}'') \in \widetilde{C} \times \widetilde{C}$.
An orbifold holomorphic map from $C$ to $X$ is represented by an orbifold holomorphic map from $\widetilde{C}$ to $X$ such that, for each node $z$,
$(\tilde{z}'$ and $\tilde{z}'')$ are mapped to $(g, x)$ and $\bigl(g^{-1},x\bigr)$, for some $g \in G$ and $x \in U$, in the twisted sector of $G \ltimes U$ which is a local model of
$\mathcal X$. The stable condition is, as usual, the finiteness of the automorphism of the map.

In ordinary Gromov--Witten theory, our focus lies on holomorphic maps from nodal Riemann surfaces.
However, in orbifold Gromov--Witten theory, we also consider holomorphic maps from
orbifold Riemann surfaces allowing orbifold nodes. It is not possible to obtain the nodal orbifold structure simply by looking at the degeneration of domain orbifold Riemann
surfaces.
Instead, we need to investigate the degeneration of holomorphic curves. To achieve this, we introduce an orbifold structure around nodes ensuring
that the homomorphism $\phi_1$ above is injective.

Similar to the case of manifolds, for a symplectic orbifold $X$, we pick an compatible almost complex structure.
Fix a homology class $A$ (the homology class of the map between coarse spaces), the moduli space
${\mathcal M}_{g, {\bf m}}(A)$ of orbifold stable maps representing $A$ is compact Hausdorff and carries a virtual fundamental class.
(The injectivity of $\phi_1$ mentioned above is used for the effectivity of the Kuranishi structure.)
Here, $g$ is the genus of the domain curve, ${{\bf m}=(m_1, \dots, m_{\ell})}$ is the data of marked orbifold points. (When $m_j=1$, the marked point is a regular point.)
We can define the evaluation map ${\rm ev}_j\colon {\mathcal M}_{g, {\bf m}}(A) \to \vert {\mathcal IX} \vert$, $j=1, \dots, \ell$.
Using them, {\it orbifold Gromov--Witten invariant}
\[{\rm GW}_{g, \ell, A}\colon \ H^*(\vert {\mathcal IX} \vert)^{\otimes \ell} \to {\mathbb Q}\]
 is defined by
\[(\alpha_1, \dots, \alpha_{\ell}) \mapsto \sum_{\vert \bf m \vert = \ell} \int_{[{\mathcal M}_{g,{\bf m}}(A)]^{\rm vir}} {\rm ev}_1^* \alpha_1 \wedge \dots \wedge {\rm ev}_{\ell}^* \alpha_{\ell}.\]
It is a convention that the grading of $H^*(\vert {\mathcal IX} \vert)$ is shifted by twice of the age of each connected component of $\vert {\mathcal IX} \vert$.

\section{Lagrangian and dihedral twisted sector}\label{sec4}
\subsection{Definition of Lagrangians}\label{sec4.1}
A Lagrangian submanifold $L$ in a symplectic manifold $X$ has a neighborhood, which is symplectomorphic to a tubular neighborhood of the zero section of the cotangent bundle
$T^*L$ of~$L$ (Weinstein). The zero section is the fixed point set of the fiberwise multiplication by $-1$, which is an anti-symplectic involution.
Hence there is a neighborhood $W$ of $L$ and an involution~$\tau\colon W \to W$ with $\tau^* \omega = - \omega$ such that $L$ is the fixed point set of $\tau$.
Based on this fact, we introduce an orientifold structure on a symplectic orbifold and define an associated Lagrangian in a symplectic orbifold as follows.

 Let $(s, t)\colon (X_1, \omega_1) \rightrightarrows ( X_0, \omega_0)$ be a proper \'etale Lie groupoid representing the symplectic orbifold $( \mc X, \omega)$, that is, $(X_0, \omega_0)$ is a symplectic manifold
and $\omega_1= s^*\omega_0 = t^*\omega_0$.

 Firstly, we introduce the notion of symplectic orientifolds.
Denote by $B\Z_2$ the action groupoid of the trivial $\Z_2 \cong \{ \pm 1 \}$ action on a point, i.e., $\{\pm 1 \} \ltimes \{ pt\}$.
 An orientifold structure on $( \mc X, \omega)$ is a proper \'etale Lie groupoid $\wt {\mc X} = \bigl(\wt X_1 \rightrightarrows \wt X_0\bigr)$
\[
 \xymatrix{ \mc X \ar[r] &\wt {\mc X} \ar[r]^\varepsilon & B\Z_2, }
 \]
where $ \wt X_0=X_0$ and $\varepsilon\colon \wt {\mc X} = \bigl(\wt X_1 \rightrightarrows \wt X_0\bigr) \to B\Z_2 \cong \{\pm 1 \} \ltimes \{ pt\}$ is a strict groupoid morphism such that $\operatorname{Ker} (\varepsilon) =\mc X$ and, for any arrow $\gamma \notin \ker  (\varepsilon)$, the local
 diffeomorphism $\psi_\gamma\colon U_{s(\gamma)} \to U_{t(\gamma)}$ is an anti-symplectomorphism. Here $U_{s(\gamma)}$ and $U_{t(\gamma)}$ are open neighbourhoods of
 $s(\gamma)$ and $t(\gamma)$ respectively in $X_0$
 such that, for $p \in U_{s(\gamma)}$ (resp.\ $q \in U_{t(\gamma)}$), $\psi_{\gamma}$ gives an arrow $\gamma'$
 (resp.\ $\gamma''$), with $s(\gamma')=p$ (resp.\ $t(\gamma'')=q$).
 Namely, $\psi_{\gamma}$ gives local sections \smash{$U_{s(\gamma)} \to \widetilde{X}_{1}$} and \smash{$U_{t(\gamma)} \to \widetilde{X}_{1}$}.
 Equivalently, we can decompose the arrow space \smash{$\wt X_1$} as a disjoint union $\wt X_1 = X_1 \cup X_{-1}$, where $X_{-1} =\bigl\{\c\in \wt X_1 \mid \e (\c) = - 1\bigr\}$.
 The induced maps $(s, t)\colon X_{-1} \rightrightarrows X_0$
 satisfy $s^*\omega_0 = - t^* \omega_0$. Note that $(s, t)\colon X_{-1} \rightrightarrows X_0$ is not a groupoid, as any product of composable arrows $\gamma$ and $\gamma'$ with $\e (\c) =\e(\c') =-1$ is not an arrow in $X_{-1}$ as $\varepsilon (\gamma\gamma') = 1$.
 We call $\gamma$ with $\varepsilon(\gamma)=-1$ an odd arrow (or odd morphism).

 For $x\in X_0$, let $\Gamma_x$ and
 $\wt \Gamma_x$ be the isotropy group of $x$ in $\mc X$ and $\wt{\mc X}$ respectively, then an orientifold structure induces an extension of the local group
 $\Gamma_x$
\[
\{1\} \to \Gamma_x \longrightarrow \wt \Gamma_x \longrightarrow  \Z_2 \to \{1 \}.
 \]
 We simply denote an orientifold
 structure on $( \mc X, \omega)$ by $\bigl( \wt{\mc X}, \omega, \e\bigr)$.

 Now, we introduce the notion of Lagrangians in symplectic orbifolds.
 Let $L $ be a subset of~$X \cong \vert \mc X \vert$.
 We shall call $L$ the underlying space of a Lagrangian in the symplectic orbifold $X$, if there is a neighborhood $W$ of $L$ equipped with
 an open orientifold structure on ${\mc W} = \pi^{-1}(W)$ in the following sense. Here $\pi\colon X_0 \to \vert \mc X \vert$ is the projection to the coarse space.

 An open orientifold structure on $(\mc X, \omega)$ is an open suborbifold $\mc W $ of $( \mc X, \omega)$ with an orientifold structure
 \smash{$\bigl( \wt {\mc W}, \omega, \e\bigr)$} on $( \mc W, \omega)$.
 Let $(\mc W, \omega)$ be an open suborbifold of $( \mc X, \omega)$, represented by~$(W_1 \rightrightarrows W_0)$, where $W_0$ is
an open submanifold of $X_0$, and $W_1 = s^{-1}(W_0) \cap t^{-1}(W_0)$. Assume that $\mc W $ admits an orientifold structure
\smash{$ \wt {\mc W} = \bigl(\bigl(\wt W_1 \rightrightarrows \wt W_0\bigr), \omega, \e\bigr)$}.
We can define a Lagrangian $\mc L = (L_1 \rightrightarrows L_0)$ in $( \mc X, \omega)$ associated to the orientifold structure \smash{$\bigl( \wt {\mc W}, \e\bigr)$} as follows,

Denote by $\operatorname{Inv}\bigl( \wt{W}_1\bigr) =\bigl\{\tau\in \wt{ W}_1 \mid \e(\tau) =-1, \, s(\tau) = t(\tau), \, \tau^2 = u( s(\tau)) \bigr\}$ the space of involutive odd arrows in \smash{$\wt{W}_1$}.
Then there is the adjoint action of $W_1$ on \smash{$\operatorname{Inv}\bigl( \wt{W}_1\bigr)$}, i.e.,
$ (\c, \tau) \in (W_1){}_s\times_{t} \operatorname{Inv}\bigl( \wt{W}_1\bigr) \mapsto \c \cdot \tau \cdot \c^{-1} \in \operatorname{Inv}\bigl( \wt{ W}_1\bigr)$.
Note that \smash{$\operatorname{Inv}\bigl( \wt{ W}_1\bigr) $} is a manifold of dimension $ \frac 12 \dim X$.
A~Lagrangian in a symplectic orbifold is the following data.
\begin{enumerate}\itemsep=0pt
\item[(1)]
$\mc I$ is a collection of connected components of \smash{$\operatorname{Inv}\bigl( \wt{ W}_1\bigr)$}, which is invariant under the adjoint action by $W_1$.
\item[(2)]
The space $L_0$ is defined by
\begin{gather*}
L_0  =  \{(x, \tau)\mid x\in W_0,\, \tau \in {\mathcal I}, \, s(\tau) = t(\tau)=x\} =  \bigsqcup_{\tau\in \mathcal I} \wt W_0^\tau,
 \end{gather*}
where \smash{$\wt W_0^\tau =\bigl\{x\in \wt{W}_0 = W_0\mid  x= s(\tau) =t(\tau)\bigr\}$} (the fixed point set of $\tau$),
such that
\begin{gather*}
L= \bigcup_{\tau\in \mathcal I} \pi \bigl(\wt W_0^\tau\bigr).
\end{gather*}

\item[(3)] The morphism space $L_1$ between two elements $(x, \tau)$ and $(y, \tau')$:
\[
{\Mor}_{\mc L} ( (x, \tau), (y, \tau'))
\]
consists of
 $\gamma \in {\Mor}_\mc W( x, y)$ satisfies $\gamma \tau \gamma^{-1} = \tau'$.
 There are canonical maps $(s, t)\colon L_1 \rightrightarrows L_0$. The composition, unit and inverse maps are induced from the corresponding maps for $\mc W =(W_1 \rightrightarrows W_0)$.
\end{enumerate}
If these conditions are fulfilled, we call $\mc L$ a {\it Lagrangian} in the symplectic orbifold $(\mc X, \omega)$ associated to an open orientifold structure \smash{$\bigl(\wt {\mc W}, \omega, \e\bigr)$} on an open suborbifold $(\mc W, \omega)$ of $(\mc X, \omega)$.
A Lagrangian in $(\mc X, \omega)$ is a Lagrangian associated to some open orientifold structure \smash{$\bigl( \wt {\mc W}, \omega, \e\bigr)$} and some $\mc I$.

Note that $\mc L =(L_1 \rightrightarrows L_0)$ is a proper \'etale Lie groupoid. There is a canonical strict morphism
\[
(\iota_1, \iota_0)\colon\ (L_1 \rightrightarrows L_0) \longrightarrow (W_1 \rightrightarrows W_0)
\]
such that $\iota_0\colon L_0 \to W_0$ is an immersion, and the groupoid structure on $ L_1 \rightrightarrows L_0$ is induced
from the groupoid action of $\mc W$ on $L_0$.

When $\mathcal W$ is locally described by a finite group action groupoid $G \ltimes U \rightrightarrows U$, the orientifold structure implies that there is an exact sequence of groups
\[ \{1\} \to G \longrightarrow \widetilde{G} \overset{\varepsilon}{\longrightarrow} \Z_2 \to \{1\} \]
such that $g \in \widetilde{G}$-action on $U$ satisfies $g^* \omega = \varepsilon(g) \omega$.
Let $I$ be a subset of \smash{$\operatorname{Inv}\bigl({\widetilde G}\bigr)$}, which is the set of all $g \in {\widetilde G}$ such that $g^2=1$ and $\varepsilon(g)=-1$.
If $I$ is invariant under the adjoint action by $G$, we can also have a Lagrangian ${\mathcal L}_I$, which is an open and closed suborbifold of
$\mathcal L$.

As an example of a Lagrangian in a symplectic orbifold, we consider
the diagonal $\Delta_\mc X$ in the product $(\mc X \times \mc X, - {\rm pr}_1^* \omega + {\rm pr}_2^* \omega)$ of a symplectic orbifold $(\mc X,\omega)$.
Firstly, we review the diagonal in the setting of orbifolds, cf.\ \cite[Example~2.6]{ALR}.
Let $\mathcal X$ be the proper \'etale Lie groupoid presenting~$X$.
Then we define the diagonal groupoid~$\Delta_{\mathcal X}$ by
\begin{align*}
(\Delta_{\mathcal X})_0  ={}&  \{(x, \alpha, y) \in X_0 \times X_1 \times X_0 \mid s(\alpha)= x,\, t(\alpha)=y \}  \\
(\Delta_{\mathcal X})_1  ={}&  \bigl\{ (x, \alpha, y) \overset{h_1, h_2}{\to} (x', \alpha', y') \mid h_1, h_2 \in X_1,  \\
 &s(h_1)=x, \, t(h_1)=x', \,s(h_2)=y,\, t(h_2)=y',\, m(h_2, \alpha) = m( \alpha', h_1) \bigr\}
\end{align*}
such that the structure maps $s$, $t$, $m$, $u$, $i$ are naturally induced from $\mathcal X$.
Then $\Delta_{\mathcal X}$ is a groupoid.
When $\mathcal X = (G \ltimes U \rightrightarrows U)$, we find that
\begin{gather*}
(\Delta_{G \ltimes U})_0 = \bigsqcup_{g \in G} \Delta_g, \Delta_g = \{(x, gx) \mid x \in U\},\qquad (\Delta_{G \ltimes U})_1 = (G \times G) \times \biggl(\bigsqcup_{g \in G} \Delta_g\biggr).
\end{gather*}
We use the notation $(x, gx; g)$ for $(x, gx) \in \Delta_g$.
We set
$s((h_1, h_2), (x, gx; g)) = (x, gx; g)$, $t((h_1,h_2), (x, gx; g))=\bigl(h_1x, h_2gx; h_2 g h_1^{-1}\bigr)$, etc.
In other words, $\Delta_{G \ltimes U}$ is the action groupoid
$(G \times G) \ltimes \sqcup \Delta_g$
associated with the action of $G \times G$ on
$\sqcup \Delta_g$.

Now we discuss the canonical orientifold structure on $(\mc X, - \omega) \times (\mc X, \omega)$ defined by the canonical
 involution $\tau_{\rm can}$ of switching two components.
 Suppose that a local model for $\mc X$ is given by a finite group action groupoid $(G\ltimes U \rightrightarrows U)$ as above.
Note that $\Delta_{\rm id}$ is the fixed point set of the involution $\tau_{\rm can}(x, x')=(x',x)$.
Since $\tau_{\rm can}$ is the reflection with respect to $\Delta_{\rm id}$, we set $\tau_{\rm id}=\tau_{\rm can}$.
Note also that $\Delta_g = (1, g) \Delta_{\rm id}$ is the fixed point set of
\[
\tau_g = (1,g) \circ \tau_{\rm id} \circ \bigl(1, g^{-1}\bigr)=\bigl(g^{-1}, g\bigr) \circ \tau_{\rm id}.
\]
Denote by \smash{$\widetilde{G\times G}$} the group generated by $G \times G$ and $\tau_{\rm id}$.
Then we see that
\smash{$\wt{G \times G}$} is the semi direct product
\[
\wt{G\times G} = (G\times G)\rtimes \Z_2
\]
with the $\Z_2$-action is defined by the adjoint action by the involution
$(g_1, g_2) \mapsto \tau_{\rm can} \cdot (g_1, g_2) \cdot \tau_{\rm can}^{-1} = (g_2, g_1)$ for $(g_1, g_2) \in G\times G.$
 Then the corresponding
local model for $\mc X\times \mc X$ is
\[
(G\times G)\ltimes ( U\times U) \rightrightarrows U \times U.
\]

Note that the odd involutive elements in $\wt {G \times G}$
consists of $\bigl\{\bigl(g^{-1}, g\bigr)\circ \tau_{\rm can} \mid g\in G\bigr\}$,
and the fixed point of the involution action $\bigl(g^{-1}, g\bigr)\circ \tau_{\rm can}$ on $ U\times U$ consists of
$
\{(x, gx) \mid x\in U\}.
$

The Lagrangian $\mc L$ for this orientifold structure can be described as follows.
\big(In this case, we take \smash{${\mathcal I} = \operatorname{Inv}\bigl( \bigl(\wt{\mc X \times \mc X}\bigr)_1\bigr)$}.\big)
For an order $m$ element $g\in G$, there is a monomorphism from the dihedral group
\[
D_m \longrightarrow \wt {G \times G}
\]
with the image generated by $\bigl(g^{-1}, g\bigr)$ and $\tau$. So, locally over $U\times U$, the associated Lagrangian is presented by $ \bigl(L_0^U \rightrightarrows L_1^U \bigr)$ with the unit space
\[L_0^U =\bigsqcup_{g\in G} \operatorname{Fix}\bigl(\bigl(g^{-1}, g\bigr)\circ \tau_{\rm can}\bigr) = \bigsqcup_{g \in G} \{ (x, gx; g) \mid g\in G,\, x\in U\} ,
\]
and the morphism space $L_1^U = (G\times G) \times L_0^U $ with the obvious action
\[
(h_1, h_2)\colon\ (x, gx; g) \mapsto \bigl(h_1 x, h_2 gx; h_2 g h_1^{-1}\bigr)
\]
 for $(h_1, h_2)\in G \times G$, here $ (h_1x, h_2 gx)$ is a fixed point of the involutive element \[\bigl(h_1g^{-1}h_2^{-1}, h_2 g h_1^{-1}\bigr)\circ \tau_{\rm can}.\]
We clearly see that $\mc L$ is equivalent to $\Delta_{\mc X}$, which we discussed above. So this canonical Lagrangian is indeed $\mc X$ diagonally embedded in $(\mc X, - \omega) \times (\mc X, \omega)$.
 Locally, we can check that the
 natural inclusion $(\phi_1, \phi_0)\colon (G \rtimes U\! \rightrightarrows U) \to \bigl(L^U_1\rightrightarrows L^U_0\bigr)$
 is the equivalence where~${\phi_0 (x) = (x, x; {\rm id})}$ and $\phi_1 (g, x) = ((g, g), (x, x; {\rm id}))$.

Clearly, the groupoid $\Delta_{G \ltimes U}$ is isomorphic to $\bigl(L^U_1 \rightrightarrows L^U_0\bigr)$, which we obtained from the symplectic orientifold.

\subsection{Dihedral twisted sector}
Let $\mc L$ be a Lagrangian of a symplectic orbifold $(\mc X, \omega)$ associated to a local orientifold structure \smash{$\bigl( \wt {\mc W}, \omega, \e\bigr)$} on an open
symplectic suborbifold $(\mc W, \omega)$ as in the previous subsection. Let $\mc W = (W_1\rightrightarrows W_0)$. We introduce a notion of the dihedral twisted sectors of a~Lagrangian in a~symplectic orbifold.

\begin{defn} For a Lagrangian $\mc L = ( L_1 \rightrightarrows L_0) $ in a symplectic orbifold $(\mc X, \omega)$ as above, we have the following proper \'etale Lie groupoid $I_{\mathcal X} {\mathcal L}$ such that
the space of objects is given by
\[(I_{\mathcal X} {\mathcal L})_0 = L_0 \times_{\wt W_0} L_0 = \bigsqcup_{\tau, \tau' \in {\mathcal I}} \wt W_0^\tau \cap \wt W_0^{\tau'} \]
and the space of morphisms is induced from the diagonal action of $\mc W$-action on $L_0 \times_{\wt W_0} L_0$ in the sense that for any point
\[
(x, \tau, \tau') \in \wt W_0^\tau \cap \wt W_0^{\tau'} \subset L_0 \times_{\wt W_0} L_0
\]
and $h\in s^{-1} (x) \subset W_1$, $h \cdot (x, \tau, \tau') = \bigl(t(h), h \tau h^{-1}, h \tau' h^{-1}\bigr).$ Equivalently, for any pair
$(x, \tau, \tau')$ and $(y, \tilde \tau, \tilde \tau')$ in $(I_{\mathcal X} {\mathcal L})_0$,
\[
{\Mor}_{I_{\mathcal X} {\mathcal L}} ((x, \tau, \tau'), (y, \tilde \tau, \tilde \tau') )=\bigl\{ h\in {\Mor}_{\mc W} (x, y) \mid
 \tilde \tau = h\tau h^{-1} ,\, \tilde \tau' = h \tau' h^{-1} \in \wt{W}_1 \bigr\}
 \]
 with the obvious source and target maps.
The other structure maps $ m, u, i$ are induced from the corresponding maps in $\mc W$.
 We call $I_{\mathcal X} {\mathcal L}$ the {\it dihedral twisted sector}\footnote{Let $\Z_m$, $D_m$ be the cyclic group of order $m$ and the dihedral group of order $2m$, respectively.
 Then the dihedral twisted sector of $\mc L$ can be described as
 \[
 \mc W\ltimes \bigsqcup_{m\in \N} \IMor \bigl(BD_m, \wt{\mc W}\bigr).\]
 Here the symbol $\IMor$ in the notation
 indicates that the injectivity on the level of morphism spaces is required.} of the Lagrangian $\mathcal L$ in the symplectic orbifold $\mathcal X$.
 \end{defn}

 Recall the exact sequence $ \{1\} \to W_1 \to \widetilde{W}_1 \overset{\varepsilon}{\to} \{\pm 1 \} \to \{1\}$. For any pair $ \tau, \tau' \in \operatorname{Inv} \bigl(\wt{\mc W}_1\bigr)$
with non-empty \smash{$ \wt W_0^\tau \cap \wt W_0^{\tau'} $}, then $\tau'\tau^{-1} = g \in W_1$ and $\tau g \tau = g^{-1}$.
Pairs $(\tau, \tau')$ of odd involutions is in one-to-one correspondence with pairs $(\tau, g)$ of odd involution and an element of $W_1$ such that $\tau \cdot g \cdot \tau = g^{-1}$.
So we can rephrase the definition of the dihedral twisted sector
using $(x, \tau, g)$ for $g\in W_1$ satisfying \smash{$\tau g\tau = g^{-1}$} (dihedral relation).\footnote{In the case of Floer theory for cleanly intersecting Lagrangians, the formulation using a pair of anti-symplectic involutions is better suited.}

Now, we discuss the case of the diagonal $\Delta_{\mathcal X}$ in Section~\ref{sec4.1}. In fact, it motivated us to make the definition of the dihedral twisted sector.
We first consider a holomorphic map $\Phi=(\Phi_0, \Phi_1)$ from a half infinite cylinder
\[
[0, \infty) \times [-1, 1]/(\tau, -1) \sim (\tau, 1)
\]
 to a complex orbifold presented by $G \ltimes U$,
the action groupoid of the action of a finite group $G$ on $U$.
Note that the half infinite cylinder is presented by the following groupoid $\mathcal C=(C_0, C_1)$ such that $C_0= [0, \infty) \times [-1, 1]$ and $C_1$ is the disjoint union of
$\{ {\rm id}_{(\tau, t)} \vert (\tau, t) \in [0, \infty) \times [-1, 1]\}$ and
$\{ a_{(\tau, 1)}\colon (\tau, 1) \to (\tau, -1)\}$, $\{a_{(\tau, -1)}\colon (\tau, -1) \to (\tau, 1)\}$. We define the structure maps $s$, $t$, $m$, $u$, $i$ in an obvious way.
This groupoid\footnote{This is not a
\'etale groupoid. Take $[0, \infty) \times (-1- \delta, 1+ \delta)$ instead of $[0, \infty) \times [-1,1]$ and an equivalence relation $(\tau, t) \sim (\tau, t + 2)$ for $t \in (-1-\delta, -1+ \delta)$.
Then the corresponding groupoid is a proper \'etale Lie groupoid representing the cylinder. For the argument here, both groupoid works.}
is a presentation of the equivalence relation $\sim$ on $C_0$.
Then~$\Phi$ is given by the pair of a holomorphic map $\Phi_0\colon [0, \infty) \times [-1, 1] \to U$, $\Phi_1\colon C_1 \to (G \ltimes U)_1$ such that $\Phi_1({\rm id}_{(\tau, t)})= u(\Phi_0(\tau, t))$,
$\Phi_1(a_{(\tau, -1)})=(g(\tau), \Phi_0(\tau, -1))$ and $\Phi_1(a_{(\tau,1)})=\bigl(g(\tau)^{-1} , \Phi_0(\tau, 1)\bigr)$ for some $g(\tau) \in G$.
Since $\Phi_1$ is continuous, $\Phi_1$ is regarded as a locally constant function with values in $G$. We set $\gamma =\Phi_1(a_{(\tau, -1)})$.
Here, the condition that $\Phi_0(\tau, 1) = \gamma \cdot \Phi_0 (\tau, -1)$ is required.
Since $\gamma$ is of finite order, we obtain a holomorphic map from a finite cover of the half-infinite cylinder to $U$.
If $\Phi_0$ is a holomorphic map with finite energy, the limit
\[
p=\lim_{\tau \to \infty} \Phi_0(\tau, t)
\]
 exists and is independent of $t$.
Thus we find that $\gamma \cdot p = p$, i.e., $(\gamma, p) \in G \ltimes U$ belongs to the twisted sector.
If we replace $\Phi$ by $\sigma \cdot \Phi=\bigl(\sigma \cdot \Phi_0, \sigma \cdot \Phi_1 \cdot \sigma^{-1}\bigr)$, $(\gamma, p)$ changes to $\bigl(\sigma \cdot \gamma \cdot \sigma^{-1}, \sigma \cdot p\bigr)$.
It is compatible with the local description of the twisted sector of the action groupoid.

Next, we rewrite the above argument using a holomorphic map
\[\Psi(\tau, t)=(\Phi_0(\tau, -t), \Phi_0(\tau, t))\]
 from
$[0, \infty] \times [0,1]$ to $(G \times G) \ltimes ((U, -J) \times (U,J))$.
Then it satisfies the boundary condition~$\Delta_{{\rm id}}$ along $t=0$ and $\Delta_{\gamma}$ along $t=1$.
Apply the Schwarz reflection principle to $\Psi$ and the anti-holomorphic involution $\tau_{\rm id} (x, x') = (x', x)$ which is the reflection with respect to $\Delta_{\rm id}$,
we obtain an extension $\Psi^+\colon [0, \infty) \times [-1,1] \to U \times U$.
Note that the anti-holomorphic involution ${\tau_{\gamma}= (1, \gamma) \circ \tau_{\rm id} \circ \bigl(1, \gamma^{-1}\bigr) = \bigl(\gamma^{-1}, \gamma\bigr) \circ \tau_{\rm id}}$,
the fixed point set of which is $\Delta_{\gamma}$.
Since $\Psi^+(\tau, 1)=\tau_{\rm id} \circ \Psi^+(\tau, -1)$, we find that
\[\Psi^+(\tau, 1) = \bigl(\gamma^{-1}, \gamma\bigr) \cdot \Psi^+(\tau, -1).\]
In this way, we get the element $\bigl(\gamma^{-1}, \gamma\bigr) \in G \times G$ from $\Psi$ and $\tau_{\rm id}$.
Then we conclude that the object $(\gamma, p)$ in the twisted sector, which is given by the behavior of $\Phi$ under $\tau \to \infty$ corresponds to
the triple $\bigl( (p,p), \tau_{\rm id}, \bigl(\gamma^{-1}, \gamma\bigr) \bigr)$ in the dihedral twisted sector, which is given by the behavior of $\Psi$ under $\tau \to \infty$.
This leads us to the notion of the dihedral twisted sector.

If we consider $\Psi_{\rho_1, \rho_2}=(\rho_1, \rho_2) \circ \Psi$, the boundary conditions along $t=0$ and $t=1$ are \smash{$\Delta_{\rho_2 \cdot \rho_1^{-1}}$} and \smash{$\Delta_{\rho_2 \cdot \gamma \cdot \rho_1^{-1}}$},
respectively. Set $\sigma= \rho_2 \cdot \rho_1^{-1}$ and $\zeta = \rho_2 \cdot \gamma \cdot \rho_1^{-1}$.
Applying the Schwarz reflection to $\Psi_{\rho_1, \rho_2}$ with respect to $\tau_{\sigma}$, we obtain an extension $\Psi_{\rho_1, \rho_2}^+ \colon [0, \infty) \times [-1, 1] \to U \times U$.
The boundary conditions for $\Psi^+_{\rho_1, \rho_2}$ are $\tau_{\sigma} \Delta_{\zeta}$ and $\Delta_{\zeta}$ along $t=-1$ and $t=1$, respectively.
Note that
\[\tau_{\zeta}=\bigl(1, \zeta \cdot \sigma^{-1}\bigr) \circ \tau_{\sigma} \circ \bigl(1, \sigma \cdot \zeta^{-1}\bigr) = \bigl(\zeta^{-1} \cdot \sigma, \zeta \cdot \sigma^{-1}\bigr) \circ \tau_{\sigma}, \]
i.e.,
\[\tau_{\rho_2 \cdot \gamma \cdot \rho_1^{-1}} = \bigl(\rho_1 \cdot \gamma^{-1} \cdot \rho_1^{-1}, \rho_2 \cdot \gamma \cdot \rho_2^{-1}\bigr) \circ \tau_{\rho_2 \cdot \rho_1^{-1}}.\]
Then, for a point $x \in \Delta_{\zeta}$, we have $x = \bigl(\zeta^{-1} \cdot \sigma, \zeta \cdot \sigma^{-1}\bigr) \circ \tau_{\sigma} (x)$.
Hence, we find that
\[\Psi_{\rho_1, \rho_2}^+ (\tau, 1) = \bigl(\zeta^{-1} \cdot \sigma, \zeta \cdot \sigma^{-1}\bigr) \Psi_{\rho_1, \rho_2}^+ (\tau, -1)
=\bigl(\rho_1\cdot\gamma^{-1} \cdot \rho_1^{-1}, \rho_2 \cdot \gamma \cdot \rho_2^{-1}\bigr) \Psi_{\rho_1, \rho_2}^+ (\tau, -1).\]
We assign an object
$\bigl((\rho_1 p, \rho_2 p), \tau_{\rho_2 \cdot \rho_1^{-1}}, \bigl(\rho_1\cdot\gamma^{-1} \cdot \rho_1^{-1}, \rho_2 \cdot \gamma \cdot \rho_2^{-1}\bigr)\bigr)$.
Note that
\[\tau_{\rho_2 \cdot \rho_1^{-1}} = \bigl(\rho_1 \cdot \rho_2^{-1}, \rho_2 \cdot \rho_1^{-1}\bigr) \circ \tau_{\rm id}=(\rho_1, \rho_2 ) \circ \tau_{\rm id} \circ \bigl(\rho_1^{-1}, \rho_2^{-1}\bigr).\]
Thus we find that the action of $(\rho_1, \rho_2) \in G \times G$ on the object space of the dihedral twisted sector sends
$\bigl((p, p), \tau_{\rm id}, \bigl(\gamma^{-1}, \gamma\bigr)\bigr)$ to \[\bigl((\rho_1 p, \rho_2 p), \tau_{\rho_2 \cdot \rho_1^{-1}}, \bigl(\rho_1\cdot\gamma^{-1} \cdot \rho_1^{-1}, \rho_2 \cdot \gamma \cdot \rho_2^{-1}\bigr)\bigr).\]
Therefore, the object in the dihedral twisted sector determined by $\Psi$ with boundary conditions $\Delta_{\rm id}$, $\Delta_{\gamma}$ and the one determined by
$(\rho_1, \rho_2) \Psi$ with boundary conditions
 \begin{align*}
 (\rho_1, \rho_2) \Delta_{\rm id} &= \Delta_{\rho_2 \cdot \rho_1^{-1}} = \operatorname{Fix} \bigl( \tau_{\rho_2 \cdot \rho_1^{-1}}\bigr),\\
(\rho_1, \rho_2) \Delta_{\gamma}& = \Delta_{\rho_2 \cdot \gamma \cdot \rho_1^{-1}}
 =\operatorname{Fix}\bigl(\bigl(\rho_1 \cdot \gamma^{-1} \cdot \rho_1^{-1}, \rho_2 \cdot \gamma \cdot \rho_2^{-1}\bigr) \circ \tau_{\rho_2 \cdot \rho_1^{-1}}\bigr)\end{align*}
 are equivalent in the sense of
Definition \ref{coarse}.

\section[Filtered A\_infty-structure associated with a Lagrangian]{Filtered $\boldsymbol{A_{\infty}}$-structure associated with a Lagrangian}\label{sec5}
In this section, we give a brief description on the filtered $A_{\infty}$-algebra associated with a Lagrangian ${\mathcal L} = (L_1 \rightrightarrows L_0)$
in a closed symplectic orbifold~$X$.
A spin structure on ${\mathcal L}$ is a spin structure on $L_0=\{(x, \tau) \in W_0 \times {\mathcal I} \mid s (\tau) = t(\tau) =x\}$, which is invariant under the action of~$W_1$.
On the dihedral twisted sector $I_{\mathcal X} {\mathcal L}$, we define a local system $\Theta$.
Recall that an object of the dihedral twisted sector is $(x, \tau, g) \in W_0 \times {\mathcal I} \times W_1$ such that $s(\tau)=t(\tau)=x$, $s(g)=t(g)=x$
and~$\tau g\tau = g^{-1}$.
This condition is equivalent to that \smash{$x\in \wt W_0^\tau \cap \wt W_0^{ g \tau}$} for \smash{$\tau, g \tau \in {\mathcal I} \subset \operatorname{Inv} \bigl(\wt{\mc W}_1\bigr)$}.
When we present $\mathcal W = (W_1 \rightrightarrows W_0)$ as a proper \'etale action Lie groupoid, we find that a connected component of $(I_{\mathcal X}{\mathcal L})_0$
is $\{x \in W_0 \mid s(\tau)=t(\tau)=x,\, s(g)=t(g)=x\} = \operatorname{Fix}(\tau) \cap \operatorname{Fix}(g \tau)$, which is a clean intersection of Lagrangian submanifolds in $W_0$.
Then we have an $O(1)$-local system $\Theta$ given in \cite[Proposition 8.1.1]{FOOO09-II}.
We find that $\Theta$ is a $O(1)$-local system on $I_{\mathcal X}{\mathcal L}$ in the sense Definition \ref{def3.3}\,(3).
Then the space $\Omega^*(I_{\mathcal X} {\mathcal L}; \Theta)$ of differential forms with coefficients in $\Theta$ is defined on the dihedral twisted sector $I_{\mathcal X} {\mathcal L}$.

The filtered $A_{\infty}$-algebra is defined on $\Omega^*(I_{\mathcal X} {\mathcal L}; \Theta) \hat{\otimes}_{\mathbb C} \Lambda_0$ by the same way as a spin Lagrangian submanifold as in Section \ref{sec2}:
\[{\mathfrak m}_k =\sum {\mathfrak m}_{k,\beta} T^{\int_\beta, \omega},\]
where ${\mathfrak m}_{0,0}=0$, ${\mathfrak m}_{1, 0} ={\rm d}$ and
\[{\mathfrak m}_{k, \beta}(\xi_1, \dots, \xi_k) = (-1)^*({\rm ev}_0)_! ({\rm ev}_1^* \xi_1 \wedge \cdots \wedge {\rm ev}_k^* \xi_k), \quad \text{unless} \quad (k,\beta)=(0, 0), (1, 0).\]
The arguments on sign and orientation bundles on the moduli spaces in \cite {On} are extended to our setting straightforwardly.
We consider pseudo-holomorphic maps from the unit disk with interior orbifold points and boundary punctures to ${\mathcal X}$
such that the boundary is mapped to $\mathcal L$ and the boundary punctures are mapped to the dihedral twisted sector in the following way.
Take a strip-like coordinate $[R, \infty) \times [0,1)$ around the puncture. Then there is $(x, \tau, g) \in (I_{\mathcal X} {\mathcal L})_0$,
the pseudo-holomorphic map is described by $w\colon [R, \infty) \times [0,1) \to W_0$ such that $w([R, \infty) \times \{0\}) \subset \operatorname{Fix}(\tau)$,
$w([R, \infty) \times \{1\}) \subset \operatorname{Fix}(g \tau)$ and $\lim_{s \to \infty} (s, t)=x$.

We adopt the argument in \cite{FOOO18,FOOOConstr2} in the setting of \cite{CDL} to construct Kuranishi structures on the moduli space of bordered orbifold stable maps.
 The notion of relative spin structures is generalized to a Lagrangian in a symplectic orbifold.
 We will have the orbifold version of Theorem~\ref{ringiso}, i.e., the orbifold quantum cohomology of a closed symplectic orbifold $X$ is isomorphic to Lagrangian Floer cohomology of the diagonal equipped with the product structure given by $\pm {\mathfrak m}_2$.
The construction of the filtered $A_{\infty}$-bimodule associated with a relative spin pair of cleanly intersecting Lagrangians is also generalized in the orbifold setting.
The details will appear elsewhere.

\subsection*{Acknowledgements}
We would like to express our gratitude to anonymous referees for their careful reading.
Bohui Chen is supported by the National Natural Science Foundation of China (no.~11890663, no.~12071322, no.~11826102, and no.~11890660), the National Key R\&D Program of China (no.~2020YFA0714000), and the Sichuan Science and Technology Program (no.~2022JDTD0019). Kaoru Ono is partly supported by JSPS Grant-in-Aid for Scientific Research 19H00636. He is also grateful to
National Center for Theoretical Sciences, Taiwan, where a part of this work is performed.

\pdfbookmark[1]{References}{ref}
\LastPageEnding

\end{document}